\documentclass[10pt]{amsart}
\usepackage{times,amsmath,amsbsy,amssymb,amscd,mathrsfs}
\usepackage{slashbox}\usepackage{graphicx,subfigure,epstopdf,wrapfig,chemarrow}

\usepackage{algorithm2e} 
\usepackage{multicol,multirow}
\usepackage{mathtools}
\usepackage[usenames,dvipsnames,svgnames,table]{xcolor}
\usepackage[all]{xy}
\usepackage{wrapfig}

\usepackage[numbered]{mcode}

\definecolor{myBlue}{rgb}{0.0,0.0,0.55}
\usepackage[pdftex,colorlinks=true,citecolor=myBlue,linkcolor=myBlue]{hyperref}

\usepackage[hyperpageref]{backref}

\usepackage{comment,enumerate,multicol,xspace}

  \newcounter{mnote}
  \let\oldmarginpar\marginpar
    \renewcommand\marginpar[1]{\-\oldmarginpar[\raggedleft\footnotesize #1]%
    {\raggedright\footnotesize #1}}

\newtheorem{theorem}{Theorem}[section]

\newtheorem{remark}[theorem]{Remark}

\newcommand{\dx}{\,{\rm d}x}
\newcommand{\dd}{\,{\rm d}}

\newcommand{\bs}{\boldsymbol}

\newcommand{\vertiii}[1]{{\left\vert\kern-0.25ex\left\vert\kern-0.25ex\left\vert #1 
    \right\vert\kern-0.25ex\right\vert\kern-0.25ex\right\vert}}

\begin{document}
\title{Programming of Finite Element Methods in MATLAB}
\author{Long Chen}

\begin{abstract}
We discuss how to implement the linear finite element method for solving the Poisson equation. We begin with the data structure to represent the triangulation and boundary conditions, introduce the sparse matrix, and then discuss the assembling process. We pay special attention to an efficient programming style using sparse matrices in MATLAB.  
\end{abstract}

\maketitle

\section{Data Structure of Triangulations}

We shall discuss the data structure to represent triangulations and boundary conditions.
\subsection{Mesh data structure} 
The matrices \mcode{node(1:N,1:d)} and \mcode{elem(1:NT,1:d+1)} are used to represent a $d$-dimensional triangulation embedded in $\mathbb  R^d$, where \mcode{N} is the number of vertices and \mcode{NT} is the number of elements. These two matrices represent two different structure of a triangulation: \mcode{elem} for the topology and \mcode{node} for the geometric embedding. 
 
The matrix \mcode{elem} represents a set of abstract simplices. The index set $\{1,2,\ldots, N\}$ is called the global index set of vertices. Here an vertex is thought as an abstract entity. For a simplex $t$, $\{1,2, \ldots, d+1\}$ is the local index set of $t$. The matrix \mcode{elem} is the mapping (pointer) from the local index to the global one, i.e., \mcode{elem(t,1:d+1)} records the global indices of $d+1$ vertices which form the abstract $d$-simplex \mcode{t}. Note that any permutation of vertices of a simplex will represent the same abstract simplex.
 
The matrix \mcode{node} gives the geometric realization of the simplicial complex. For example, for a 2-D triangulation, \mcode{node(k,1:2)} contains $x$- and $y$-coordinates of the $k$-th nodes, respectively. 

The geometric realization introduces an ordering of the simplex. For each \mcode{elem(t,:)}, we shall always order the vertices of a simplex such that the signed area is positive. That is in 2-D, three vertices of a triangle is ordered counter-clockwise and in 3-D, the ordering of vertices follows the right-hand rule. 

\begin{remark}\rm
Even with the orientation requirement, certain permutation of vertices is still allowed. Similarly any labeling of simplices in the triangulation, i.e. any permutation of rows of \mcode{elem} matrix will represent the same triangulation. The ordering of simplexes and vertices will be used to facilitate the implementation of the local mesh refinement and coarsening. See \mcode{bisect}, \mcode{coarsen}, \mcode{bisect3}, and \mcode{coarsen3} in $i$FEM~\cite{Chen:2009Integrated}. $\qed$
\end{remark}

\begin{figure}[htp]
\begin{center}
\subfigure[A triangulation of a L-shape domain.]{
\begin{minipage}[t]{0.45\linewidth}
\centering
\includegraphics[width=4.5cm]{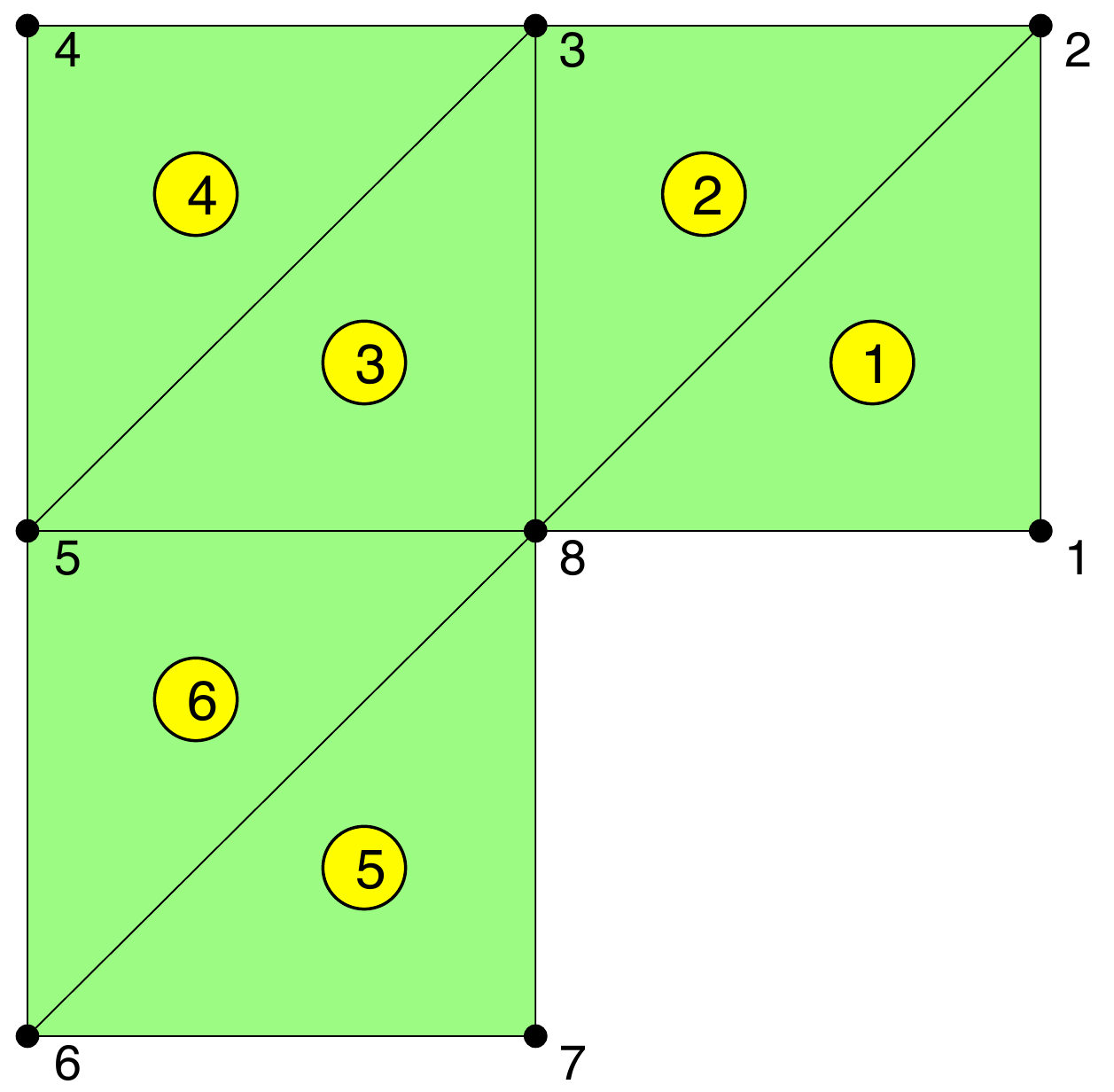}
\end{minipage}}
\subfigure[\mcode{node} and \mcode{elem} matrices]
{\begin{minipage}[t]{0.55\linewidth}
\centering
\includegraphics*[width=2.5in]{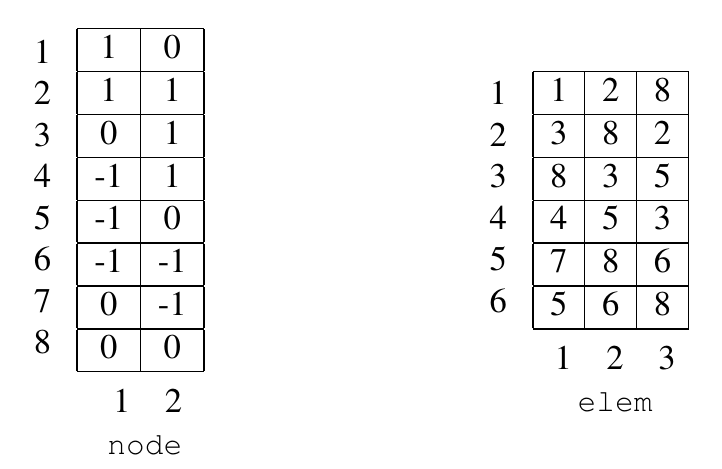}
\end{minipage}}
\caption{(a) A triangulation of the L-shape domain $(-1,1)\times (-1,1) \backslash ([0,1]\times [0,-1])$. (b) \mcode{node} and \mcode{elem} matrices.}
\label{fig:Lshape}
\end{center}
\end{figure}

As an example, \mcode{node} and \mcode{elem} matrices for a triangulation of the L-shape domain $(-1,1)\times (-1,1) \backslash ([0,1]\times [0,-1])$ are given in the Figure \ref{fig:Lshape} (a) and (b). 

\subsection{Boundary conditions}\label{sec:bd}
We use \mcode{bdFlag(1:NT,1:d+1)} to record the type of boundary sides (edges in 2-D and faces in 3-D). The value is the type of boundary condition: 
\begin{itemize}
\item 0 for non-boundary sides;
\item 1 for the first type, i.e., Dirichlet boundary;
\item 2 for the second type, i.e., Neumann boundary;
\item 3 for the third type, i.e., Robin boundary.
\end{itemize}
For a $d$-simplex, we label its $(d-1)$-faces in the way so that the $i$th face is opposite to the $i$th vertex. For example, for a 2-D triangulation, \mcode{bdFlag(t,:) = [1 0 2]} means, the edge opposite to \mcode{elem(t,1)} is a Dirichlet boundary edge, the one to \mcode{elem(t,3)} is of Neumann type, and the other is an interior edge. 

We may extract boundary edges for a 2-D triangulation from \mcode{bdFlag} by:
\begin{lstlisting}
totalEdge = [elem(:,[2,3]); elem(:,[3,1]); elem(:,[1,2])];
Dirichlet = totalEdge(bdFlag(:) == 1,:);
Neumann = totalEdge(bdFlag(:) == 2,:); 
\end{lstlisting}


\begin{remark}\rm
The matrix \mcode{bdFlag} is sparse but we use a dense matrix to store it. It would save storage if we record boundary edges or faces only. The current form is convenient for the local refinement and coarsening since the boundary conditions can be easily updated along with the change of elements. 
We do not save \mcode{bdFlag} as a sparse matrix since updating a sparse matrix is time consuming. We set up the type of \mcode{bdFlag} to \mcode{int8} to minimize the waste of spaces. $\qed$
\end{remark}

\section{Sparse matrix in MATLAB}

MATLAB is an interactive environment and high-level programming language for numeric scientific computation. One of its distinguishing features is that the main data type is the matrix. Matrices may be manipulated element-by-element, as in low-level languages like Fortran or C. But it is better to manipulate matrices at a time which will be called {\it high level} coding style. This style will result in more compact code and usually improve the efficiency.

We start with explanation of the sparse matrix and corresponding operations. The fast sparse matrix package and built-in functions in MATLAB will be used extensively later on. The content presented here is mostly based on Gilbert, Moler and Schereiber~\cite{Gilbert.J;Moler.C;Schreiber.R1992}.

One of the nice features of finite element methods is the sparsity of the matrix obtained via the discretization. Although the matrix is $N\times N=N^2$, there are only $c N$ nonzero entries in the matrix with a small constant $c$. Sparse matrix is the corresponding data structure to take advantage of this sparsity. Sparse matrix algorithms require less computational time by avoiding operations on zero entries and sparse matrix data structures require less memory by not storing many zero entries. We refer to the book~\cite{Pissanetsky.S1984} for detailed description on sparse matrix data structure and~\cite{Saad.Y2003} for a quick introduction on popular data structures of sparse matrix. In particular, the sparse matrix data structure and operations has been added to MATLAB by Gilbert, Moler and Schereiber and documented in~\cite{Gilbert.J;Moler.C;Schreiber.R1992}. 

\subsection{Storage schemes}
There are different types of data structures for the sparse matrix. All of them share the same basic idea: use a single array to store all nonzero entries and two additional integer arrays to store the indices of nonzero entries.

An intuitive scheme, known as {\it coordinate format}, is to store both the row and column indices. In the sequel, we suppose $A$ is a $m\times n$ matrix containing only $nnz$ nonzero elements. Let us look at the following simple example:
\begin{equation}\label{sparseA}
A =
\left [ \; 
\begin{matrix} 
1 & 0 & 0\, \\ 
0 & 2  & 4\, \\ 
0 & 0 & 0\,\\
0 & 9  & 0
\end{matrix} 
\right ],
\quad
i = 
\left [ \; 
\begin{matrix} 
1 \, \\ 
2 \, \\ 
4 \, \\ 
2 \, 
\end{matrix} 
\right ],
\quad
j = 
\left [ \; 
\begin{matrix} 
1 \, \\ 
2 \, \\ 
2 \, \\ 
3 \, 
\end{matrix} 
\right ],
\quad
s =
\left [ \; 
\begin{matrix} 
1 \, \\ 
2 \, \\ 
9 \, \\ 
4 \, 
\end{matrix} 
\right ].
\end{equation}
In this example, $i$ vector stores row indices of non-zeros, $j$ column indices, and $s$ the value of non-zeros. All three vectors have the same length $nnz$. The two indices vectors $i$ and $j$ contains redundant information. We can compress the column index vector $j$ to a column pointer vector with length $n+1$. The value $j(k)$ is the pointer to the beginning of $k$-th column in the vector of $i$ and $s$, and $j(n+1)=nnz+1$. For example, in CSC formate, the vector to store the column pointer will be $j=[\, 1 \; 2 \; 4 \; 5\, ]^{\intercal}$. This scheme is known as {\it Compressed Sparse Column (CSC)} scheme and is used in MATLAB sparse matrices package. 

Comparing with coordinate formate, CSC formate saves storage for $nnz - n -1$ integers which could be non-negligilble when the number of nonzero is much larger than that of the column. In CSC formate it is efficient to extract a column of a sparse matrix. For example, the $k$-th column of a sparse matrix can be build from the index vector $i$ and the value vector $s$ ranging from $j(k)$ to $j(k+1)-1$. There is no need of searching index arrays. An algorithm that builds up a sparse matrix one column at a time can be also implemented efficiently~\cite{Gilbert.J;Moler.C;Schreiber.R1992}. On the other hand, extract one row of a sparse matrix saved in CSC will involve a search in the whole vector $i$ and $j$.

\begin{remark} \rm
CSC is the internal representation of sparse matrices in MATLAB. For the convenience of users, the coordinate scheme is presented as the interface. This allows users to create and decompose sparse matrices in a more straightforward way. $\Box$
\end{remark}

Comparing with the dense matrix, the sparse matrix lost the direct connection between the index \mcode{(i,j)} and the physical location to store the value \mcode{A(i,j)}. Then accessing and manipulating one element at a time requires the searching of the index vectors to find such nonzero entry. It takes time at least proportional to the logarithm of the length of the column; inserting or removing a nonzero may require extensive data movement~\cite{Gilbert.J;Moler.C;Schreiber.R1992}. Therefore, {\it do not manipulate or change a sparse matrix element-by-element in a large \mcode{for} loop in MATLAB.}

Due to the lost of the link between the index and the value of entries, the operations on sparse matrices is delicate. One needs to code specific subroutines for standard matrix operations: matrix times vector, addition of two sparse matrices, and transpose of sparse matrices etc. Since some operations will change the sparse pattern, typically there is a priori loop to set up the nonzero pattern of the resulting sparse matrix. Good sparse matrix algorithms should follow the rule ``time is proportional to flops"~\cite{Gilbert.J;Moler.C;Schreiber.R1992}: The time required for a sparse matrix operation should be proportional to the number of arithmetic operations on nonzero quantities. 
 
\subsection{Create and decompose sparse matrix}
To create a sparse matrix, we first form $i,j$ and $s$ vectors, i.e., a list of nonzero entries and their indices, and then call the function \mcode{sparse} using $i,j,s$ as input. Several alternative forms of \mcode{sparse} (with more than one argument) allow this. The most commonly used one is

\mcode{A = sparse(i,j,s,m,n)}.

This call generates an $m\times n$ sparse matrix, having one nonzero for each entry in the vectors $i,j$, and $s$. The first three arguments all have the same length. However, the indices in $i$ and $j$ need not be given in any particular order and could have duplications. If a pair of indices occurs more than once in $i$ and $j$, \mcode{sparse} adds the corresponding values of $s$ together. This nice summation property is very useful for the assembling procedure in finite element computation. 

The function \mcode{[i,j,s]=find(A)} is the inverse of \mcode{sparse} function. It will extract the nonzero elements together with their indices. The indices set $(i,j)$ are sorted in column major order and thus the nonzero \mcode{A(i,j)} is sorted in lexicographic order of \mcode{(j,i)} not \mcode{(i,j)}. See the example in (\ref{sparseA}).

\begin{remark} \rm
There is a similar command  \mcode{accumarray} to  create a dense matrix $A$ from indices and values. It usage is slightly different from \mcode{sparse}. The index \mcode{[i j]} should be paired together to form a subscript vectors. So is the dimension \mcode{[m n]}. Since the accessing of a single element in a dense matrix is much faster than that in a sparse matrix, when $m$ or $n$ is small, say $n=1$, it is better to use \mcode{accumarray} instead of \mcode{sparse}. A most commonly used command is

\mcode{accumarray([i j], s, [m n])}.
\end{remark}

\section{Assembling of The Matrix Equation}\label{ch:assembling}
In this section, we discuss how to obtain the matrix equation for the linear finite element method for solving the Poisson equation
\begin{equation}\label{PoissonEq}
- \Delta u = f  \hbox{ in }  \Omega, \qquad u  = g_D  \hbox{ on } \Gamma _D,  \qquad  \nabla u\cdot n = g_N  \hbox{ on } \Gamma _N,
\end{equation}
where $\partial \Omega = \Gamma _D\cup \Gamma _N$ and $\Gamma _D\cap \Gamma _N=\emptyset$. We assume $\Gamma _D$ is closed and $\Gamma _N$ open.

Denoted by $H_{g,D}^1(\Omega)=\{v\in L^2(\Omega), \nabla v\in L^2(\Omega) \hbox{ and } v|_{\Gamma _D} = g_D\}$. Using integration by parts, the weak form of the Poisson equation (\ref{PoissonEq}) is: find $u\in H_{g,D}^1(\Omega)$ such that
\begin{equation}\label{weakform}
a(u,v) := \int _{\Omega} \nabla u\cdot \nabla v \dd \bs x = \int _{\Omega} fv \, \dd \bs x + \int _{\Gamma _N} g_N v \,{dS} \quad \text{ for all } v\in H_{0,D}^1(\Omega).
\end{equation}

Let $\mathcal T$ be a triangulation of $\Omega$. We define the linear finite element space on $\mathcal T$ as
$$
\mathbb  V_{\mathcal T} = \{v\in C(\bar \Omega) : v|_{\tau}\in \mathcal P_1, \forall \tau \in \mathcal T\},
$$
where $\mathcal P_1$ is the space of linear polynomials. For each vertex $v_i$ of $\mathcal T$, let $\phi _i$ be the piecewise linear function such that $\phi _i(v_i)=1$ and $\phi _i(v_j)=0$ if $j\neq i$. Then it is easy to see $\mathbb  V_{\mathcal T}$ is spanned by $\{\phi _i\}_{i=1}^{N}$. The linear finite element method for solving (\ref{PoissonEq}) is to find $u\in \mathbb  V_{\mathcal T}\cap H_{g,D}^1(\Omega)$  such that (\ref{weakform}) holds for all $v\in \mathbb  V_{\mathcal T}\cap H_{0,D}^1(\Omega)$. 

We shall discuss an efficient way to obtain the algebraic equation. It is an improved version, for the sake of efficiency, of that in the paper~\cite{Alberty.J;Carstensen.C;Funken.S1999}. 

\subsection{Assembling the stiffness matrix}
For a function $v\in \mathbb  V_{\mathcal T}$, there is a unique representation: $v=\sum _{i=1}^Nv_i\phi _i$.  We define an isomorphism $\mathbb  V_{\mathcal T}\cong \mathbb R^N$ by
\begin{equation}\label{vecrep}
v=\sum _{i=1}^Nv_i\phi _i \longleftrightarrow \bs v=(v_1, \cdots, v_N)^{\intercal},
\end{equation}
and call $\bs v$ the coordinate vector of $v$ relative to the basis $\{\phi _i\}_{i=1}^{N}$. Following the terminology in the linear elasticity, we introduce the {\it stiffness matrix}
$$
\bs A=(a_{ij})_{N\times N}, \, \hbox{ with } \quad a_{ij}=a(\phi _j,\phi _i).
$$
In this subsection, we discuss how to form the matrix $\bs A$ efficiently in MATLAB.

\subsubsection{Standard assembling process}
By the definition, for $1\leq i,j\leq N$,
$$
a_{ij} = \int _{\Omega} \nabla \phi _j\cdot \nabla \phi _i \, {\rm d\bs x}= \sum _{\tau \in \mathcal T} \int _{\tau}\nabla \phi _j\cdot \nabla \phi _i\, {\rm d\bs x}.
$$
For each simplex $\tau$, we define the local stiffness matrix $A^{\tau} = (a_{ij}^{\tau})_{(d+1)\times (d+1)}$ as 
$$
a_{i_{\tau}j_{\tau}}^{\tau} = \int _{\tau}\nabla \lambda_{j_\tau}\cdot \nabla \lambda _{i_{\tau}} \, {\rm d\bs x}, \hbox{ for } 1\leq i_{\tau}, j_{\tau}\leq d+1.
$$
The computation of $a_{ij}$ will then be decomposed into the computation of local stiffness matrix and the summation over all elements. Here we use the fact that restricted to one simplex, the basis $\phi_i$ is identical to the barycentric coordinate $\lambda_{i_{\tau}}$ and the subscript in $a_{i_{\tau}j_{\tau}}^{\tau}$ is the local index while in $a_{ij}$ it is the global index.
The assembling process is to distribute the quantity associated to the local index to that to the global index.

Suppose we have a subroutine \mcode{locatstiffness} to compute the local stiffness matrix, to get the global stiffness matrix, we apply a \mcode{for} loop of all elements and distribute element-wise quantity to node-wise quantity. A straightforward MATLAB code is like

\begin{lstlisting}
function A = assemblingstandard(node,elem)
N=size(node,1); NT=size(elem,1);
A=zeros(N,N); %A = sparse(N,N);
for t=1:NT
	At=locatstiffness(node(elem(t,:),:));
	for i=1:3
		for j=1:3
			A(elem(t,i),elem(t,j))=A(elem(t,i),elem(t,j))+At(i,j);
		end
	end
end
\end{lstlisting}

The above code is correct but not efficient. There are at least two reasons for the slow performance.
\begin{enumerate}
\item The stiffness matrix \mcode{A} created in line 3 is a dense matrix which needs $\mathcal O(N^2)$ storage. It will be out of memory quickly when $N$ is big  (e.g., $N=10^4$). Sparse matrix should be used for the sake of memory. Nothing wrong with MATLAB. Coding in other languages also need to use the sparse matrix data structure. Use \mcode{A = sparse(N,N)} in line 3 will solve this problem.

\item There is a large \mcode{for} loops with size \mcode{NT} the number of elements. This can quickly add significant overhead when \mcode{NT} is large since each line in the loop will be interpreted in each iteration. This is the weak point of MATLAB. Vectorization should be applied for the sake of efficiency.
\end{enumerate}

We now discuss the standard procedure: transfer the computation to a reference simplex through an affine map, on computing of the local stiffness matrix. We include the two dimensional case here for the comparison and completeness.

We call the triangle $\hat \tau$  spanned by $\hat{v}_1=(1,0), \hat{v}_2=(0,1)$ and $\hat{v}_3=(0,0)$ {\it a reference triangle} and use $\hat{\bs x}=(\hat x, \hat y)^{\intercal}$ for the vector in the reference coordinate. For any $\tau \in \mathcal T$, we treat it as the image of $\hat \tau$ under an affine map: $F: \hat \tau \to \tau$. One of such affine map is to match the local indices of three vertices, i.e., $F(\hat v_i)=v_i, i=1,2,3$:
$$
F(\hat{\bs x}) = B^{\intercal} (\hat{\bs x}) + c,
$$
where
$$
B = \left [
\begin{matrix}
x_1 - x_3 & y_1 - y_3 \\
x_2 - x_3 & y_2 - y_3
\end{matrix}
\right ],
\; \hbox{ and }\;
c = (x_3,  y_3)^{\intercal}.
$$
We define $\hat u(\hat{\bs x}) = u(F(\hat{\bs x}))$. Then $\hat \nabla \hat u = B\nabla u$ and ${\rm dxdy} = |\det(B)| {\rm d\hat xd\hat y}$. By change of variables, the integral becomes
\begin{align*}
\int _{\tau} \nabla \lambda _i \cdot \nabla \lambda _j {\rm dxdy}
& =
\int _{\hat \tau} (B^{-1}\hat \nabla \hat{\lambda}_i) \cdot (B^{-1}\hat \nabla \hat{\lambda}_j ) |\det(B)| {\rm d\hat xd\hat y}\\
& =
\frac{1}{2}|\det(B)|(B^{-1}\hat \nabla \hat{\lambda}_i) \cdot (B^{-1}\hat \nabla \hat{\lambda}_j ).
\end{align*}
In the reference triangle, $\hat{\lambda} _1= \hat x, \hat{\lambda} _2 = \hat y$ and $\hat{\lambda} _3= 1 - \hat x - \hat y$. Thus
$$
\hat \nabla \hat{\lambda} _1 =
\left [
\begin{matrix}
1 \\
0
\end{matrix}
\right ],
\;
\hat \nabla \hat{\lambda} _2 =
\left [
\begin{matrix}
0 \\
1
\end{matrix}
\right ],
\;
\hbox{ and }
\hat \nabla \hat{\lambda} _3 =
\left [
\begin{matrix}
-1 \\
-1
\end{matrix}
\right ].
$$
We then end with the following subroutine~\cite{Alberty.J;Carstensen.C;Funken.S1999} to compute the local stiffness matrix in one triangle $\tau$.
\begin{lstlisting}
function [At,area] = localstiffness(p)
At = zeros(3,3); 
B = [p(1,:)-p(3,:); p(2,:)-p(3,:)]; 
G = [[1,0]',[0,1]',[-1,-1]'];
area = 0.5*abs(det(B));
for i = 1:3
   for j = 1:3
      At(i,j) = area*((B\G(:,i))'*(B\G(:,j)));
   end
end
\end{lstlisting}
The advantage of this approach is that by modifying the subroutine \mcode{localstiffness}, one can easily adapt to new elements and new equations. 

\subsubsection{Assembling using sparse matrix}

As we mentioned before, updating one single element of a sparse matrix in a large loop is very expensive since the nonzero indices and values vectors will be reformed and a large of data movement is involved. Therefore the code in line 8 of \mcode{assemblingstandard} will dominate the whole computation procedure. 
In this example, numerical experiments show that the subroutine \mcode{assemblingstandard} behaves like an $\mathcal O(N^2)$ algorithm.

We should call \mcode{sparse} command once to form the sparse matrix. The following subroutine is suggested by T. Davis~\cite{Davis.T}.

\begin{lstlisting}
function A = assemblingsparse(node,elem)
N = size(node,1); NT = size(elem,1);
i = zeros(9*NT,1); j = zeros(9*NT,1); s = zeros(9*NT,1);
index = 0;
for t = 1:NT
    At = localstiffness(node(elem(t,:),:));
    for ti = 1:3
        for tj = 1:3
            index = index + 1;
            i(index) = elem(t,ti);
            j(index) = elem(t,tj);
            s(index) = At(ti,tj);
        end
    end
end
A = sparse(i, j, s, N, N);
\end{lstlisting}

In the subroutine \mcode{assemblingsparse}, we first record a list of index and nonzero entries inside the loop and use built-in function \mcode{sparse} to form the sparse matrix outside of the loop. By doing in this way, we avoid updating a sparse matrix inside a large loop. The subroutine \mcode{assemblingsparse} is much faster than \mcode{assemblingstandard}. 
This simple modification is recommended when translating C or Fortran codes into MATLAB.

\subsubsection{Vectorization}
There is still a large loop in the subroutine \mcode{aseemblingsparse}. We shall use the vectorization technique to avoid the outer large \mcode{for} loop. 

Given a $d$-simplex $\tau$, recall that the barycentric coordinates $\lambda _j(\bs x),j=1,\cdots, d+1$ are linear functions of $\bs x$. If the $j$-th vertex of a simplex $\tau$ is the $k$-th vertex of the triangulation, then the hat basis function $\phi _k$ restricted to a simplex $\tau$ will coincide with the barycentric coordinate $\lambda _j$. Note that the index $j=1,\cdots, d+1$ is the local index set for the vertices of $\tau$, while $k=1,\cdots, N$ is the global index set of all vertices in the triangulation. 

We shall derive a formula for $\nabla \lambda _i, i=1,\cdots, d+1$. Let $F_i$ denote the $(d-1)$-face of $\tau$ opposite to the $i$th-vertex. Since $\lambda _i(\bs x)=0$ for all $\bs x\in F_i$, and $\lambda _i(\bs x)$ is an affine function of $\bs x$, the gradient $\nabla \lambda _i$ is a normal vector of the face $F_i$ with magnitude $1/h_i$, where $h_i$ is height, i.e., the distance from the vertex $x_i$ to the face $F_i$. Using the relation $|\tau| = \frac{1}{d} |F_i|h_i$, we end with the following formula
\begin{equation}\label{dlambda}
\nabla \lambda _i = \frac{1}{d!\, |\tau|}\bs n_i,
\end{equation}
where $\bs n_i$ is an {\it inward} normal vector of the face $F_i$ with magnitude $\|\bs n_i\|=(d-1)!|F_i|$. Therefore
$$
a_{ij}^{\tau} = \int _{\tau} \nabla \lambda _i \cdot \nabla \lambda _j\, \dd \bs x = \frac{1}{d!^2|\tau|}\bs n_i\cdot \bs n_j.
$$

Note that we do not normalize the normal vector since the scaled one is easier to compute. In 2-D, the scaled normal vector $\bs n_i$ can be easily computed by a rotation of the edge vector.  For a triangle spanned by $\bs x_1,\bs x_2$ and $\bs x_3$, we define $\bs l_i = \bs x_{i+1} - \bs x_{i-1}$ where the subscript is 3-cyclic. For a vector $\bs v = (x,y)$, we denoted by $\bs v^{\bot} = (-y, x)$. Then $\bs n_i = \bs l_i ^{\bot}$ and $\bs n_i\cdot \bs n_j=\bs l_i \cdot \bs l_j$. The edge vector $\bs l_i$ for all triangles can be computed using matrix operations and can be used to calculate the area of all triangles. 

We end with the following compact, efficient, and readable code for the assembling of stiffness matrix in two dimensions.

\begin{lstlisting}
function A = assembling(node,elem)
N = size(node,1); NT = size(elem,1);
ii = zeros(9*NT,1); jj = zeros(9*NT,1); sA = zeros(9*NT,1);
ve(:,:,3) = node(elem(:,2),:)-node(elem(:,1),:);
ve(:,:,1) = node(elem(:,3),:)-node(elem(:,2),:);
ve(:,:,2) = node(elem(:,1),:)-node(elem(:,3),:);
area = 0.5*abs(-ve(:,1,3).*ve(:,2,2)+ve(:,2,3).*ve(:,1,2));
index = 0;
for i = 1:3
    for j = 1:3
        ii(index+1:index+NT) = elem(:,i); 
        jj(index+1:index+NT) = elem(:,j);
        sA(index+1:index+NT) = dot(ve(:,:,i),ve(:,:,j),2)./(4*area);
        index = index + NT;
    end
end
A = sparse(ii,jj,sA,N,N);
\end{lstlisting}

\begin{remark}\rm
One can further improve the efficiency by using the symmetry of the matrix. For example, the inner loop can be changed to \mcode{for j = i:3}. $\Box$ 
\end{remark}

In 3-D, the scaled normal vector $\bs n_i$ can be computed by the cross product of two edge vectors. We list the code below and explain it briefly.

\begin{lstlisting}
function A = assembling3(node,elem)
N = size(node,1); NT = size(elem,1);
ii = zeros(16*NT,1); jj = zeros(16*NT,1); sA = zeros(16*NT,1);
face = [elem(:,[2 4 3]);elem(:,[1 3 4]);elem(:, [1 4 2]);elem(:, [1 2 3])];
v12 = node(face(:,2),:)-node(face(:,1),:);
v13 = node(face(:,3),:)-node(face(:,1),:);
allNormal = cross(v12,v13,2);
normal(1:NT,:,4) = allNormal(3*NT+1:4*NT,:);
normal(1:NT,:,1) = allNormal(1:NT,:);
normal(1:NT,:,2) = allNormal(NT+1:2*NT,:);
normal(1:NT,:,3) = allNormal(2*NT+1:3*NT,:);
v12 = v12(3*NT+1:4*NT,:); 
v13 = v13(3*NT+1:4*NT,:);
v14 = node(elem(:,4),:)-node(elem(:,1),:);
volume = dot(cross(v12,v13,2),v14,2)/6;
index = 0;
for i = 1:4
    for j = 1:4
        ii(index+1:index+NT) = elem(:,i); 
        jj(index+1:index+NT) = elem(:,j);
        sA(index+1:index+NT) = dot(normal(:,:,i),normal(:,:,j),2)./(36*volume);
        index = index + NT;
    end
end
A = sparse(ii,jj,sA,N,N);
\end{lstlisting}
The code in line 4 will collect all faces of the tetrahedron mesh. So the \mcode{face} is of dimension \mcode{4NT$\times$3}. For each face, we form two edge vectors \mcode{v12} and \mcode{v13}, and apply the cross product to obtain the scaled normal vector in \mcode{allNormal} matrix. The code in line 8-11 is to reshape the \mcode{4NT$\times$3} normal vector to a \mcode{NT$\times$3$\times$4} matrix. Note that in line 8, we assign the value to \mcode{normal(:,:,4)} first such that the MATLAB will allocate enough memory for the array \mcode{normal} when creating it. Line 15 use the mix product of three edge vectors to compute the volume and line 19--22 is similar to 2-D case. The introduction of the scaled normal vector $\bs n_i$ simplify the implementation and enable us to vectorize the code.

\subsection{Right hand side}
We define the vector $\bs f = (f_1,\cdots, f_N)^{\intercal}$ by $f_i=\int _{\Omega}f\phi _i$, where $\phi _i$ is the hat basis at the vertex $v_i$. For quasi-uniform meshes, all simplices are around the same size, while in adaptive finite element method, some elements with large mesh size could remain unchanged. Therefore, although the 1-point quadrature is adequate for the linear element on quasi-uniform meshes, to reduce the error introduced by the numerical quadrature, we compute the load term $\int _{\Omega}f \phi _i$ by 3-points quadrature rule in 2-D and 4-points rule in 3-D. Arbitrary order quadrature will be discussed in the next section.

We list the 2-D code below to emphasize that the command \mcode{accumarray} is used to avoid the slow \mcode{for} loop over all elements.
\begin{lstlisting}
mid1 = (node(elem(:,2),:)+node(elem(:,3),:))/2;
mid2 = (node(elem(:,3),:)+node(elem(:,1),:))/2;
mid3 = (node(elem(:,1),:)+node(elem(:,2),:))/2;
bt1 = area.*(f(mid2)+f(mid3))/6;
bt2 = area.*(f(mid3)+f(mid1))/6;
bt3 = area.*(f(mid1)+f(mid2))/6;
b = accumarray(elem(:),[bt1;bt2;bt3],[N 1]);
\end{lstlisting}

\subsection{Boundary conditions}\label{sec:bd}
We list the code for 2-D case and briefly explain it for the completeness. Recall that \mcode{Dirichlet} and \mcode{Neumann} are boundary edges which can be found using \mcode{bdFlag}. See Section 1.
\begin{lstlisting}
%-------------------- Dirichlet boundary conditions------------------------
isBdNode = false(N,1); 
isBdNode(Dirichlet) = true;
bdNode = find(isBdNode);
freeNode = find(~isBdNode);
u = zeros(N,1); 
u(bdNode) = g_D(node(bdNode,:));
b = b - A*u;
%-------------------- Neumann boundary conditions -------------------------
if (~isempty(Neumann))
    Nve = node(Neumann(:,1),:) - node(Neumann(:,2),:);
    edgeLength = sqrt(sum(Nve.^2,2)); 
    mid = (node(Neumann(:,1),:) + node(Neumann(:,2),:))/2;
    b = b + accumarray([Neumann(:),ones(2*size(Neumann,1),1)], ... 
                   repmat(edgeLength.*g_N(mid)/2,2,1),[N,1]); 
end
\end{lstlisting}

Line 2-4 will find all Dirichlet boundary nodes. The Dirichlet boundary condition is posed by assign the function values at Dirichlet boundary nodes \mcode{bdNode}. It could be found by using \mcode{bdNode = unique(Dirichlet)} but \mcode{unique} is very costly. So we use logic array to find all nodes on the Dirichlet boundary, denoted by \mcode{bdNode}. The other nodes will be denoted by \mcode{freeNode}.

The vector \mcode{u} is initialized as zero vector. Therefore after line 7, the vector \mcode{u} will represent a function $u_D\in H_{g, D}$. Writing $u=\tilde u + u_D$, the problem (\ref{weakform}) is equivalent to finding $\tilde u\in \mathbb  V_{\mathcal T}\cap H_0^1(\Omega)$ such that $a(\tilde u, v) = (f, v) -a(u_D,v)+ (g_N,v)_{\Gamma _N}$ for all $v\in \mathbb  V_{\mathcal T}\cap H_0^1(\Omega)$. The modification of the right hand side $ (f, v) -a(u_D,v)$ is realized by the code \mcode{b=b-A*u} in line 8. The boundary integral involving the Neumann boundary edges is computed in line 11--15 using the middle point quadrature. Note that it is vectorized using \mcode{accumarray}.
 
Since $u_D$ and $\tilde u$ use disjoint nodes set, one vector \mcode{u} is used to represent both. The addition of $\tilde u + u_D$ is realized by assign values to different index sets of the same vector \mcode{u}. We have assigned the value to boundary nodes in line 5. We will compute $\tilde u$, i.e., the value at \mcode{freeNode}, by the direct solver
\begin{equation}\label{directsolve}
\mcode{u(freeNode)=A(freeNode,freeNode)$\backslash$ b(freeNode)}.
\end{equation}

For the Poisson equation with pure Neumann boundary condition
$$
- \Delta u = f  \; \text{ in } \Omega,  \quad \frac{\partial u}{\partial n}=g \text{ on }\Gamma,
$$
there are two issues on the well posedness of the continuous problem:
\begin{enumerate}
 \item solutions are not unique. If $u$ is a solution of Neumann problem, so is $u+c$ for any constant $c\in \mathbb  R$. One more constraint is needed to determine this constant. A common choice is $\int_{\Omega} u \dx = 0$.
 
 \item a compatible condition for the existence of a solution. There is a compatible condition for $f$ and $g$:
\begin{equation}\label{compatible}
-\int _{\Omega}f \, dx = \int _{\Omega}\Delta u \dx = \int _{\partial \Omega} \frac{\partial u}{\partial n} \dd S = \int _{\partial \Omega} g \dd S.
\end{equation}
\end{enumerate}

We discuss the consequence of these two issues in the discretization. 
For Neumann problem, the stiffness matrix \mcode{A} is symmetric but only semi-definite. The kernel of \mcode{A} consists of constant vectors, i.e, the rank of \mcode{A} is \mcode{N-1}. Then \mcode{Au=b} is solvable if and only if 
\begin{equation}\label{meanb}
\mcode{mean(b)=0} 
\end{equation}
which is the discrete compatible condition. If the integral is computed exactly, according to \eqref{compatible}, \eqref{meanb} should hold in the discrete case. Since numerical quadrature is used to approximate the integral, \eqref{meanb} may hold exactly. We can enforce \eqref{meanb} by the modification \mcode{b = b - mean(b)}.

To deal with the constant kernel of \mcode{A}, we can simply set \mcode{freeNode=2:N} and then use \eqref{directsolve} to find values of \mcode{u} at \mcode{freeNode}. Since solution \mcode{u} is unique up to a constant, afterwards we can modify \mcode{u} to satisfy certain constraint. For example, to impose the zero average, i.e., $\int_{\Omega}u\dx = 0$, we could use the following code:
\begin{lstlisting}
c = sum(mean(u(elem),2).*area)/sum(area);
u = u - c;
\end{lstlisting}
The $H^1$ error will not affect by the constant shift but when computing $L^2$ error, make sure the exact solution will satisfy the same constraint. 
\section{Numerical Quadrature} 

In the implementation, we need to compute various integrals on a simplex.  In this section, we will present several numerical quadrature rules for
simplexes in 1, 2 and 3 dimensions. 

The numerical quadrature is to approximate an integral by weighted average of function values at sampling points $p_i$:
$$
\int_{\tau} f(\bs x)\dd \bs x \approx I_n(f):=\sum _{i=1}^{n} f(p_i)w_i |\tau|.
$$
The order of a numerical quadrature is defined as the largest integer $k$ such that $\int_{\tau} f = I_n(f)$ when $f$ is a polynomial of degree less than equal to $k$. 

A numerical quadrature is determined by the quadrature points and corresponding weight: $(p_i, w_i), i=1,\ldots, n$. For a $d$-simplex $\tau$, let $\bs x_i, i=1,\ldots, d+1$ be vertices of $\tau$. The simplest one is the one point rule: 
$$
I_1(f) = f(c_{\tau})|\tau|, \quad c_{\tau} = \frac{1}{d+1}\sum_{i=1}^{d+1}\bs x_i.
$$
A very popular one is the trapezoidal rule: 
$$
I_1(f) = \frac{1}{d+1}\sum_{i=1}^{d+1}f(\bs x_i)|\tau|.
$$
Both of them are of order one, i.e., exact for the linear polynomial only. For second order quadrature, in 1-D, the Simpson rule is quite popular
$$
\int_{a}^bf(x)\dx \approx (b-a)\frac{1}{6}\left ( f(a)+ 4f((a+b)/2)+f(b)\right ).
$$
For a triangle, a second order quadrature, i.e., exact for quadratic polynomials, is using three middle points $m_i, i=1,2,3$ of edges:
$$
\int_{\tau} f(\bs x) \dd\bs x \approx \frac{|\tau|}{3}\sum_{i=1}^3 f(m_i).
$$
These rules are popular due to the reason that the points and the weight are easy. No such second order rule exists for a tetrahedron in 3-D.

A criterion for choosing quadrature points is to attain a given precision with the fewest possible function evaluations. For the two (center v.s. vertices) first order quadrature rules given above, which one is more efficient? Restricting to one element, the answer is the center. When considering the evaluation over the whole triangulation, the trapezoidal rule is better since it only evaluates the function at $N$ vertices while the center rule needs $NT$ evaluation. It is a simple exercise to show $NT \approx 2N$ asymptotically in 2-D. 


In 1-D, the Gauss quadrature use $n$ points to achieve the order $2n-1$ which is the highest order for $n$ points. The Gauss points are roots of orthogonal polynomials and can be found in almost all text books on numerical analysis.   
Quadrature rules for triangles and tetrahedron which is less well documented in the literature and we refer to~\cite{Zhang.L;Cui.T;Liu.H2009} for a set of symmetric quadrature rules. $16$ digits accurate quadrature points are included in $i$FEM. Type \mcode{quadpts} and \mcode{quadpts3} for the usage. We present the points in the barycentric coordinate $p = (\lambda_1, \ldots, \lambda_{d+1})$. The Cartesian coordinate of $p$ is obtained by $\sum _{i=1}^{d+1}\lambda_i\bs x_i$. 

\bibliographystyle{abbrv}

\end{document}